\documentclass[10pt,a4paper]{article}
\usepackage[T1]{fontenc}
\usepackage[utf8]{inputenc}
\usepackage{lmodern,amsmath,amssymb,booktabs,array,microtype}
\usepackage[margin=19.5mm]{geometry}
\usepackage[hidelinks]{hyperref}
\usepackage{titlesec}
\titlespacing*{\section}{0pt}{1.2ex plus .4ex}{.55ex}
\hypersetup{pdftitle={Thirteen Improved Binary Constant-Weight Codes with Explicit Cyclic Orbit Certificates},pdfauthor={Peter Rafat Adel Iskandar}}
\newcommand{\A}{\mathrm{A}}
\begin{document}
\begin{center}
{\Large\bfseries Thirteen Improved Binary Constant-Weight Codes\\[2pt] with Explicit Cyclic Orbit Certificates}\par
\vspace{.45em}
Peter Rafat Adel Iskandar\\[-.15em]
{\small Independent researcher, Egypt \quad \href{mailto:1peteradel@gmail.com}{1peteradel@gmail.com}}
\end{center}
\vspace{-.6em}
\begin{abstract}
We exhibit thirteen explicit binary constant-weight codes improving lower bounds recorded in Brouwer's table, including
$\A(29,6,5)\ge333$, $\A(55,4,5)\ge66{,}144$, and $\A(62,4,5)\ge107{,}360$.
Their 408,729 codewords are specified by 11,386 orbit representatives under explicitly given cyclic coordinate permutations, together with a standalone reconstruction and distance verifier.
All thirteen orbit expansions were checked to agree exactly with the full codeword lists.
\end{abstract}

\section{Results}
Let $\A(n,d,w)$ be the largest cardinality of a binary length-$n$, weight-$w$ code with pairwise Hamming distance at least $d$.
Tables of bounds on $\A(n,d,w)$ go back to \cite{Brouwer1990} and are maintained online by Brouwer \cite{BrouwerTable}.
Recent computational improvements include \cite{Rosin2026,Echols2026,Dzoba2026,Lysen2026}. In particular, \cite{Dzoba2026} constructs codes through prescribed coordinate automorphisms and orbit-based search, an approach related to the orbit certificates presented here.
We provide thirteen improved lower bounds, recorded on 9 November 2026 in Brouwer's online table under the label \texttt{Ad} (Peter Adel) \cite{BrouwerTable}.
Table~\ref{tab:codes} gives the previous recorded lower bounds, the new certified sizes, and an explicit cyclic presentation of every code.
The column $U$ reports valid upper bounds from Johnson recursion; for $(29,6,5)$ we use the sharper tabulated bound 364. These upper bounds are not asserted to be the best available.

\begin{table}[htbp]
\centering\small
\setlength{\tabcolsep}{4.0pt}
\caption{Explicit orbit presentations and certified lower bounds. $M_0$ is the table value recorded immediately before the 9 November submission, $M$ is the new size, $g$ is described by its consecutive cycle lengths, $r$ is the number of orbit representatives, and $U$ is the reported upper bound. All rows satisfy $M\le \A(n,d,w)\le U$.}
\label{tab:codes}
\begin{tabular}{@{}rrr r r c r r@{}}
\toprule
$n$&$d$&$w$&$M_0$&$M$&cycle lengths of $g$&$r$&$U$\\\midrule
16&4&7&836&\textbf{841}&\(7^{2},1^{2}\)&121&1,040\\
20&4&5&874&\textbf{876}&\(3^{6},1^{2}\)&294&912\\
25&4&6&7,787&\textbf{7,800}&\(5^{5}\)&1,560&8,379\\
29&4&5&4,466&\textbf{4,491}&\(3^{9},1^{2}\)&1,503&4,750\\
29&6&5&315&\textbf{333}&\(3^{9},1^{2}\)&111&364\\
36&4&5&11,032&\textbf{11,207}&\(35,1\)&321&11,397\\
38&4&5&13,566&\textbf{13,737}&\(19^{2}\)&723&14,341\\
39&4&5&15,276&\textbf{15,618}&\(38,1\)&411&16,450\\
44&4&5&25,564&\textbf{25,585}&\(43,1\)&595&26,488\\
54&4&5&59,943&\textbf{60,632}&\(53,1\)&1,144&61,959\\
55&4&5&65,164&\textbf{66,144}&\(53,1^{2}\)&1,248&68,013\\
60&4&5&93,206&\textbf{94,105}&\(59,1\)&1,595&95,748\\
62&4&5&106,445&\textbf{107,360}&\(61,1\)&1,760&109,678\\
\bottomrule
\end{tabular}
\end{table}

\section{Cyclic construction and certificates}
For each row of Table~\ref{tab:codes}, label the coordinates $0,\ldots,n-1$. The cycle-length column partitions them into \emph{consecutive blocks}.
Define $g$ to cyclically shift every block of length $q>1$ forward by one position, simultaneously, fixing singleton blocks. For instance, the $(55,4,5)$ construction uses
\[
 g=(0\ 1\ \cdots\ 52)(53)(54),\qquad |\langle g\rangle|=53.
\]
The accompanying \texttt{ORBIT\_MANIFEST.csv} gives the blocks and the corresponding file of one binary representative per selected orbit. We define the corresponding code by
\[
\mathcal C_{n,d,w}=\bigcup_{s\in R_{n,d,w}}\{g^j(s):0\le j<|\langle g\rangle|\}.
\]
The orbits in each union are disjoint. Their cardinalities, reported in the manifest, sum to $M$. The full automorphism groups were not determined; each $\langle g\rangle$ is a verified subgroup of the automorphism group of its code.
For example, the $(55,4,5)$ construction is exactly the union of 1,248 disjoint orbits of length 53, hence has 66,144 words; the $(36,4,5)$ construction has 320 orbits of length 35 and one of length 7, hence has 11,207 words.
Every generated word has weight $w$. For even $d$, no $(w-d/2+1)$-element subset of coordinates is contained in two distinct generated supports; therefore each code has minimum distance at least $d$.
The standard-library script \texttt{verify\_orbit\_codes.py} regenerates each full code from the generators and representatives, checks its exact size and minimum-distance condition, and can optionally confirm set equality against the full codeword lists.

\section{Bounds and provenance}
For any constant-weight code, the Johnson recurrence gives
\[
\A(n,d,w)\le \left\lfloor\frac{n}{w}\A(n-1,d,w-1)\right\rfloor.
\]
For $d=4$ and $w=5$, Table~\ref{tab:codes} uses the established exact values of $\A(n-1,4,4)$ tabulated by Brouwer \cite{BrouwerTable}; for $w=6,7$ the recurrence is iterated from exact smaller-weight values. At $(29,6,5)$ we use the sharper upper bound $364$ recorded in Brouwer's table. For example, the exact value $\A(54,4,4)=6183$ \cite{BrouwerTable} gives
\[
66{,}144\le \A(55,4,5)\le \left\lfloor\frac{55}{5}\cdot 6183\right\rfloor=68{,}013.
\]

The codes were obtained by computational search. The cyclic presentations in Table~\ref{tab:codes} were subsequently \emph{recovered and verified from the final word lists}; they are provided as concise, reproducible certificates of the resulting constructions.

The ancillary files supply the cyclic generators and 11,386 orbit representatives in \texttt{ORBIT\_MANIFEST.csv}, checksums, a reconstruction verifier, and complete lists of the 408,729 codewords. The supplementary material verifies the resulting codes but does not reproduce the original search procedure.

\textbf{Disclosure.} Anthropic Claude Code and OpenAI ChatGPT assisted with computational experimentation, certificate software, structural analysis, literature checking, and manuscript preparation under the author's direction. The author takes responsibility for the results and their release.

\end{document}